\documentclass{amsart}
\usepackage[dvips]{graphicx}
\usepackage{pgfplots}
\usepackage{amssymb}
\usepackage{verbatim}
\usepackage{graphicx}
\usepackage{amsmath}
\usepackage{kpfonts}
\usepackage{tikz}

\newcommand{\Beq}{\begin{equation}}
\newcommand{\Eeq}{\end{equation}}

 \newtheorem{theorem}{Theorem}[section]
\newtheorem{lemma}[theorem]{Lemma}
\newtheorem{prop}[theorem]{Proposition}
\newtheorem{exam}[theorem]{Example}
\newtheorem{corollary}[theorem]{Corollary}
\newtheorem{definition}[theorem]{Definition}

\begin{document}


\title{Turbulent closed relations}
\author{Judy Kennedy, Christopher Mouron, Van Nall}
\maketitle  

\begin{abstract}
This paper generalizes the classical notion of turbulence from dynamical systems generated by continuous functions to those defined by closed relations on compact metric spaces. Using the Mahavier product and the associated shift map, we introduce and explore CR-turbulence and reverse CR-turbulence, analyzing their relationship to topological entropy. A key focus is understanding when turbulence implies entropy and vice versa, with results showing that for finite closed relations, these properties are equivalent. However, examples are provided to demonstrate that this equivalence can fail for more general relations.

We also construct a large class of explicit turbulent closed relations on the unit interval that are dynamically rich yet structurally simple. Additionally, since homeomorphisms cannot admit turbulence, we investigate a weakened notion of turbulence: separated, continuum-wise semi-turbulence for homeomorphisms, and then prove that smooth fans cannot support even this weaker form of turbulent dynamics. The paper includes new examples, counterexamples, and open questions that deepen the understanding of turbulence in non-classical settings.

\end{abstract}
\-

\section{Introduction} 

Turbulence was first introduced in the well-written and accessible text \textbf{Dynamics in One Dimension}, by L. Block and W. A. Coppel \cite{BlockCoppel}, in order to study the dynamics of maps of the interval to itself. They state there that the terminology was suggested by a paper by A. Lasota and J.A. Yorke \cite{LasotaYorke}. Block and Coppel first proved that turbulent maps on the interval have periodic points of all periods. They define a map $f$ of a closed interval to itself to be ``chaotic'' if $f^n$  is separated turbulent for some $n$. (Since the work ``chaos'' now has so many meanings in mathematics, we stick with the term ``turbulence'' and variations of that idea, here.) They connected the notion to other ``chaotic'' dynamical properties - in particular, to entropy and periodic points. 

\vphantom{}

This paper was motivated by a lovely talk given by Veronica Mart\'inez de la Vega at the Summer Topology Conference in Coimbra, Portugal, in 2024. The results of her recent research on turbulence with co-author Alejandro Illanes can be found in \cite{illanesVeronica}.

\vphantom{}

 The study of the dynamics of a continuous function on an interval is particularly well developed, and the book of Block and Coppel does a wonderful job of summarizing the development of this study at the time the book was published in 1992. This study is not only of interest in itself, it is also necessary for an understanding of dynamical systems on higher-dimensional spaces. 
 
 \vphantom{}
 
Here we generalize the notion of turbulence to closed relation dynamical systems and study the properties of such systems. Each closed relation $F$ on a compact metric space $X$ gives rise to a subspace $X_F^+$ of $X^{\infty}$. The natural shift map $\sigma: X^{\infty} \to X^{\infty}$ is invariant on the subspace $X_F^+$, i.e., $\sigma(X_F^+) \subset X_F^+$. The behavior of $\sigma|X_F^+$ is intimately related to the closed relation on $F$, and gives us a tool with which to study the dynamics of a closed relation on $X$. In particular, we study the connection between positive entropy on closed relations and various kinds of turbulence on closed relations. 

\vphantom{}

The paper is organized as follows:

\vphantom{}

Section 2 lays the foundational concepts, notation, and terminology used throughout the paper. It defines turbulence (and variants) in classical and generalized contexts, introduces closed relations and Mahavier products, and provides background on related topics such as  fans and topological entropy.

\vphantom{}

Section 3 investigates the relationship between different notions of turbulence for closed relations and shift maps on Mahavier products. The goal is to bridge classical turbulence concepts (for functions) and the new CR-turbulence (closed relation turbulence) defined for closed relations. 

\vphantom{}

Section 4 focuses on closed relations $F \subset X \times X$ where F is a finite set, and it investigates the conditions under which such relations are CR-turbulent and/or have positive entropy. We show that, for finite closed relations, a variety of important dynamical properties are equivalent.

\vphantom{}

Sections 3 and 4 establish that for many closed relations (especially finite ones), positive entropy and CR-turbulence are closely related or even equivalent. Section 5 challenges that idea by showing that positive entropy does not always imply CR-turbulence or reverse CR-turbulence. In other words, entropy and turbulence can diverge in more general (non-finite) closed relations.

\vphantom{}

From the results of M. Misiurewicz (see \cite{Misiurewicz}), it is known that for continuous maps of the interval, the map has positive topological entropy if and only if some iterate of the map is separated turbulent. There is a complete dsicussion of these and related results in Chapter VIII of \cite{BlockCoppel}. So, a natural question arises: If a closed relation $F \subset [0,1]^2$ is the union of graphs of two continuous functions defined on $[0,1]$, does $ent(F) > 0$ imply that $F$ is CR-turbulent or reverse CR-turbulent? We present a counterexample (see Example 5.3) if one of the function domains is not connected. 

\vphantom{}

Section 6 presents explicit examples of closed relations on the unit interval [0,1] that are CR-turbulent, reverse CR-turbulent, and structurally simple, yet dynamically rich. The goal is to provide a constructive and general method for generating such turbulent relations.

\vphantom{}

It is not possible for homeomorphisms to have turbulence, but they can possess a related weaker property. Section 7 investigates this weaker turbulence property in homeomorphisms on a class of topological spaces called smooth fans. We call this weaker form of turbulence separated, continuum-wise semi-turbulence, and investigate whether such turbulence is possible for homeomorphisms on smooth fans.

\section{Definitions, Notation and Background}

\begin{definition} Let $X$ be a compact metric space. We use $d$ or $\rho$ to denote the metric on $X$, unless otherwise specified. We use $\mathbb{N}$ to denote the positive integers; $\mathbb{N}_0$ to denote the nonnegative integers; and $\mathbb{Z}$ to denote the integers.

\end{definition}

\vphantom{}

\begin{definition} A \emph{continuum} is a nonempty compact connected metric space. A continuum is \emph{degenerate} if it consists of exactly one point. Otherwise it is \emph{nondegenerate}. A subspace of a continuum which is itself a continuum is a \emph{subcontinuum}. 
	
\end{definition}

\vphantom{}

\begin{definition} Let $X$ and $Y$ be metric spaces, and let $f:X \to Y$ be a function. We use $\Gamma(F)= \{(x,y): y=f(x) \}$ to denote the \emph{graph of the function} $f$.
	
\end{definition}

\vphantom{}

 \begin{definition}
 Let $X$ be a compact metric space and let $G \subset X \times X$ be a relation on $X$. If $G$ is closed in $X \times X$, then $G$ is a \emph{closed relation on} $X$.  
 \end{definition}

\vphantom{}

\begin{definition} Let $X$ be a set and let $G$ be a relation on $X$. Then we define $G^{-1}=\{ (y,x): (x,y) \in G \}$. Also, $G^{-1}$ is the \emph{inverse relation of the relation $G$ on} $X$.
	
\end{definition}

\vphantom{}

\begin{definition} Let $X$ be a compact metric space and let $G$ be a closed relation on $X$. 

\begin{itemize}
\item  We define the \emph{Mahavier product} $G \star 
G$ as  $\{(x,y,z): (x,y) \in G, (y,z) \in G\}$. If $\mathbf{a}=(x,y) \in G$  and $\mathbf{b}=(y,z) \in G$, then we define $\mathbf{a} \star \mathbf{b}$ to be the point $\mathbf{a} \star \mathbf{b}=(x,y,z) \in G \star G$.

\item If $m$ is a positive integer, we define the $mth$ \emph{Mahavier product} $\star_{i=1}^m G$ to be 
$\star_{i=1}^m G = \{(x_1,x_2, \ldots, x_{m+1}) \in \Pi_{i=1}^{m+1} X: \text{ for each } i \in \{1,2,\ldots,m\} \}, (x_i,x_{i+1}) \in G \}.$ We abbreviate $\star_{i=1}^m G$ by $X_G^m$.
\item We define the  \emph{infinite Mahavier product} $\star_{i=1}^{\infty} G$ to be \newline 
$\star_{i=1}^{\infty} G = \{(x_1,x_2, \ldots) \in \Pi_{i=1}^{\infty} X: \text{ for each positive integer } i, (x_i,x_{i+1}) \in G \}.$ We abbreviate $\star_{i=1}^{\infty} G$ by $X_G^+$.

\item We define the function $$\sigma_G^+: X_G^+ \to X_G^+$$ given by $$\sigma_G^+(x_1,x_2, \ldots) = (x_2,x_3, \ldots)$$ for each 
$(x_1,x_2, \ldots) \in X_G^+$. The map $\sigma_G^+$ is called the \emph{shift map on} $X_G^+$, and it is a continuous function. (To see this, observe that $X_G^+$ is invariant under the shift $\sigma$ on $\Pi_{i=1}^{\infty} X$ defined by $\sigma(x_1,x_2, \dots)=(x_2,x_3, \dots)$ for $(x_1,x_2,\dots) \in \Pi_{i=1}^{\infty} X$. Since $\sigma$ is continuous, so is $\sigma|X_G^+ = \sigma_G^+$.)
\end{itemize}

\end{definition}

\vphantom{}

\begin{definition}(\emph{The standard projection} $\pi_i$.) Suppose $X$ is a space. If $n$ is a positive integer, and $\mathbf{x}=(x_1,x_2, \ldots, x_n) \in X^n$, then $\pi_i(\mathbf{x}) = x_i$ for $i \in \{1,2,\dots,n\}$. If $\mathbf{x}=(x_1,x_2, \ldots) \in X^{\infty}$, then $\pi_i(\mathbf{x}) = x_i$ for $i$ each positive integer $i$. If $m \le n$, $m$ and $n$ are integers, we use $[m,n]$ to denote the set $\{m,m+1,\ldots,n\}$ and refer to $[m,n]$ as an integer interval (if no confusion arises). Also, if $\mathbf{x}=(x_1,x_2,\dots)$, then $\pi_{[m,n]}(\mathbf{x})=(x_m,x_{m+1},\ldots, x_n)$, when $0 <m \le n$, $m,n$ integers.
	
\end{definition}

\vphantom{}

\begin{definition} If $K \subset X^n$ for some positive integer $n$, then the set $K^{-1}$ is defined to be the set $$K^{-1} = \{(x_n,x_{n-1}, \ldots, x_1): (x_1, x_2, \dots, x_n) \in K \}.$$ We call $K^{-1}$ the \emph{inverse} of the set $K$.
	
\end{definition}

\vphantom{}

\begin{definition} Suppose $X$ is a compact metric space and $f:X \to X$ is continuous. Then $(X,f)$, or just $f$ when there is no confusion, is a \emph{dynamical system}.
	
\end{definition}

\vphantom{}

\begin{definition} A dynamical system $(X,f)$ is \emph{turbulent} if there are nonempty closed sets $K$ and $L$ in $X$ such that $K \cap L$ contains at most one point and $K \cup L \subset f(K) \cap f(L)$.
	
\end{definition}

\vphantom{}

\begin{definition} A dynamical system $(X,f)$ is \emph{separated turbulent} if there are disjoint nonempty closed sets $K$ and $L$ in $X$ such that $K \cup L \subset f(K) \cap f(L)$.
	
\end{definition}

\vphantom{}

\begin{definition} Suppose $n$ is a positive integer. A dynamical system $(X,f)$ is  \emph{ $n$-turbulent}  if  
$(X,f^n)$ is turbulent.

\end{definition}

\vphantom{}

\begin{definition} A dynamical system $(X,f)$ is  \emph{ separated $n$-turbulent}  if  
$(X,f^n)$ is separated turbulent.
\end{definition}

\vphantom{}

\begin{prop} Suppose $X$ is a space and $f:X \to X$ is a function. Disjoint nonempty closed sets $K_0$ and $K_1$ in $X$ have the property that $K_0 \cup K_1 \subset f(K_0) \cap f(K_1)$ if and only if for each sequence $i_0,i_1, \ldots$ of $0$'s and $1$'s, $K_{i_0} \subset f(K_{i_1}) \subset f^2(K_{i_2}) \subset \cdots$.

\end{prop}

\vphantom{}

\begin{definition} A closed relation $F \subset X \times X$ is \emph{CR-turbulent} if there are $n \in \mathbb{N}$ and two disjoint nonempty closed sets $K$ and $L$ in $X_{F}^{n}$  such that $\pi_1(K) \cup \pi_1(L) \subset \pi_n(K) \cap \pi_n(L)$.
	\end{definition}

\vphantom{}

\begin{definition} A closed relation $F \subset X \times X$ is \emph{reverse CR-turbulent} if there are $n \in \mathbb{N}$ and two disjoint nonempty closed sets $K$ and $L$  in $X_F^{n}$ such that $\pi_n(K) \cup \pi_n(L) \subset \pi_1(K) \cap \pi_1(L)$.

\end{definition}

\vphantom{}

Since fans and topological entropy appear later in the paper, we give a brief, basic discussion of these terms here. More details appear later, as they are needed.

\subsection{Fans} 

At first glance, fans seem deceptively simple. Their structure - arcs radiating from a single point - makes them feel intuitively easy to understand. However, the formal definition of a fan is surprisingly intricate, and this complexity is not accidental. It reflects the need to precisely capture subtle topological properties that arise in these spaces.

\vphantom{}

To illustrate this, consider the following open question:
Is every compact metric space X satisfying the conditions below necessarily a fan?

\begin{enumerate}
\item	$X$ is one-dimensional;
\item	There exists a point $v \in X$ and a family $\mathcal{F}$ of arcs in $X$ such that:

\begin{enumerate}
\item For any distinct $A, B \in \mathcal{F}$, we have $A \cap B = \{v\}$;
\item $X = \bigcup \{ F : F \in \mathcal{F} \}$.

\end{enumerate}

\end{enumerate}

\vphantom{}

Intuitively, such a space looks like a fan - a union of arcs joined only at a common point - but whether it satisfies the formal definition remains an open problem. 
   
\begin{definition}
Let $X$ be a continuum. 
\begin{enumerate}
\item The continuum $X$ is \emph{unicoherent} if for any subcontinua $A$ and $B$ of $X$ such that $X=A\cup B$,  the compactum $A\cap B$ is connected. 
\item The continuum $X$ is \emph{hereditarily unicoherent } provided that each of its subcontinua is unicoherent.
\item The continuum $X$ is a \emph{dendroid} if it is an arcwise connected, hereditarily unicoherent continuum
\item  Let $X$ be a dendroid.  A point $x\in X$ is called an \emph{end-point of $X$} if for  every arc $A$ in $X$ that contains $x$, $x$ is an end-point of $A$.  The set of all end-points of $X$ is denoted by $E(X)$.  
\item The continuum $X$ is \emph{a simple triod} if it is homeomorphic to $([-1,1]\times \{0\})\cup (\{0\}\times [0,1])$.
\item Let $X$ be a simple triod. A point $x\in X$ is called \emph{the top-point} or, briefly, the \emph{top of $X$} if  there is a homeomorphism $\varphi:([-1,1]\times \{0\})\cup (\{0\}\times [0,1])\rightarrow X$ such that $\varphi(0,0)=x$.
\item Let $X$ be a dendroid.  A point $x\in X$ is called \emph{a ramification-point of $X$}, if there is a simple triod $T$ in $X$ with top $x$.  The set of all ramification-points of $X$ is denoted by $R(X)$. 
\item The continuum $X$ is a \emph{fan} if it is a dendroid with at most one ramification point $v$, which is called the \emph{top} of the fan $X$ (if it exists).
\item  Let $X$ be a fan.  We say that $X$ is \emph{a Cantor fan} if $X$ is homeomorphic to the continuum $\bigcup_{c\in C}A_c$, where $C\subseteq [0,1]$ is the standard Cantor set and for each $c\in C$, $A_c$ is the  {convex} segment in the plane from $(0,0)$ to $(c,1)$.
\item A fan $X$ is a \emph{smooth} fan if it can be embedded in the Cantor fan.
\item A \emph{Lelek fan} is a smooth fan that has a dense set of endpoints. It is known that Lelek fans are unique (see \cite{oversteegen} and \cite{jjc1}), i.e., any two Lelek fans are homeomorphic.
\end{enumerate}	
\end{definition}

\subsection{Topological entropy}

Topological entropy, or simply entropy as used here, is a fundamental but often nonintuitive concept in dynamical systems. For a continuous function from a compact metric space to itself, topological entropy is either a non-negative number or $\infty$. Roughly speaking, it quantifies the rate at which distinct orbits of the function become increasingly distinguishable over time - essentially, how much the function “mixes up” points in the space. Though defined globally, entropy often reflects highly local dynamical complexity, making it a subtle but powerful invariant.

\vphantom{}

Topological entropy has been extensively studied across many areas of mathematics as a rigorous measure of ``chaotic'' behavior. In this paper, we adopt the definition originally introduced by Adler, Konheim, and McAndrew, presented clearly in Chapter 7 of Peter Walters’ book \cite{Walters book}. While several equivalent formulations exist, they all yield the same numerical value for a given system, assuming the same logarithm base is used.

 \begin{definition}
Let $X$ be a set,  let $f:X\rightarrow X$ be a function and let $\mathcal S$ be a family of subsets of $X$. Then we define
$$
f^{-1}(\mathcal S)=\{f^{-1}(S) \ | \ S\in \mathcal S\}.
$$
\end{definition}

\begin{definition}
Let $X$ be a set and let $\mathcal A_1$, $\mathcal A_2$, $\mathcal A_3$, $\ldots$, $\mathcal A_m$ be families of subsets of $X$.  Then we define
$$
\vee_{i=1}^{m}\mathcal A_i=\{A_1\cap A_2\cap A_3\cap \ldots \cap A_m \ | \ \textup{ for each } i\in \{1,2,3,\ldots ,m\}, A_i\in \mathcal A_i\}.
$$
\end{definition}
\begin{definition}
Let $X$ be a compact metric space and let $f:X\rightarrow X$ be a continuous function. For any open cover $\alpha$ for $X$, we define
$$
ent(f,\alpha)=\lim_{m\to\infty}\frac{\log(N(\vee_{i=0}^{m}f^{-i}(\alpha)))}{m}.
$$
 
\end{definition}

\begin{definition}
Let $X$ be a compact metric space and let $f:X\rightarrow X$ be a continuous function. We define the \emph{topological entropy} $ent(f)$ to be
$$
ent(f) = \mbox{  } \sup_{\alpha} \medspace ent(f,\alpha)
$$
where $\alpha$ ranges over all open covers of $X$.
	
\end{definition}

 \section{Turbulence and CR-turbulence}
In this section, we address whether every turbulent dynamical system is also separated  $n$-turbulent for some $n \in \mathbb{N}$, and explore the relationship between a closed relation $F$ on a compact metric space $X$ with CR-turbulence and the various turbulence properties of the shift map $\sigma_F^+$ on the infinite Mahavier product $X_F^+$. We begin the section by giving a result about the relationship of CR-turbulence for a closed relation and separated turbulence for a shift map on the Mahavier product determined by that relation.

\vphantom{}

\begin{theorem}

A closed relation of $F \subset X \times X$ is CR-turbulent if and only if there is an $n \in \mathbb{N}$, $n \ge 2$, such that $(X_F, (\sigma_F^+)^{n-1})$ is separated turbulent.

\end{theorem} 

\begin{proof}  

Assume the relation $F\subset X \times X$ is CR-turbulent. Then there is an $n \in \mathbb{N}$ and there are disjoint nonempty closed sets $K_0$ and $K_1$ in $X^n_F \subset X^{n+1}$ such that $\pi_1(K_0) \cup \pi_1(K_1) \subset \pi_n(K_0) \cap \pi_n(K_1)$. For each $i \in \mathbb{N}$, let $r(i)=[(i-1)(n-1)+1,i(n-1)+1]$.

\vphantom{}

We use $\Sigma_2$ to denote the collection of all sequences $\alpha= \alpha_1, \alpha_2, \ldots$ such that $\alpha_n \in \{0,1\}$ for each positive integer $n$. Let $\alpha \in \Sigma_2$ and for each $j \in \mathbb{N}$, define 
$$ M_{j,\alpha} = \{ \mathbf{x} \in X^{\infty}: \text{ for each } 1\le i \le j, \pi_{r(i)}(\mathbf{x}) \in K_{\alpha(i)} \}. $$

We show that for each $j \in \mathbb{N}$ and each $\alpha \in \Sigma_2$, $M_{j,\alpha}$ is a nonempty closed set. The proof is by induction on $j$: If $j=1$, then $M_{1,\alpha}= \{\mathbf{x} \in X^{\infty}:\pi_{[1,n]}(\mathbf{x}) \in K_{\alpha(1)} \}$, which is closed and nonempty whether $\alpha(1) =0$ or $\alpha(1) =1$. Now assume $l \in \mathbb{N}$ and $l >1$, and $M_{l-1,\alpha} \neq \emptyset$ for all $\alpha \in \Sigma_2$. Choose $\alpha \in \Sigma_2$ and let $\alpha ' \in \Sigma_2$ such that $\alpha '(i)=\alpha(i+1)$ for each $i \in \mathbb{N}$. Let $\mathbf{x} \in M_{l-1,\alpha'}$.  Since $\pi_1(\mathbf{x}) \in \pi_1(K_0) \cup \pi_1(K_1) \subset \pi_n(K_0) \cap \pi_n(K_1)$, there is a $\mathbf{y} \in K_{\alpha(1)}$ such that $\pi_n(\mathbf{y})=\pi_1(\mathbf{x})$. It follows that $\mathbf{y} \star \mathbf{x} \in M_{l,\alpha}$. So $M_{j,\alpha} \neq \emptyset$ for each $j \in \mathbb{N}$.

\vphantom{}

Clearly, $M_{j+1,\alpha} \subset M_{j,\alpha}$ for each $j \in \mathbb{N}$. So, since $X^{\infty}$ is compact, $\cap_{j=1}^{\infty} M_{j,\alpha} \neq \emptyset$. Let $M_{\alpha} = \cap_{j=1}^{\infty} M_{j,\alpha}$. Note that $M_{\alpha} \subset X_F^+$ for each $\alpha \in \Sigma_2$. Let $M=\cup_{\alpha \in \Sigma_2} M_{\alpha}$, let $K_0^*=\{\mathbf{x} \in M: \pi_{[1,n]} \in K_0 \}$ and let $K_1^*=\{\mathbf{x} \in M: \pi_{[1,n]} \in K_1 \}$. Then $K_0^*$ and $K_1^*$ are nonempty closed sets in $X_F^+$ such that $(\sigma_F^+)^{n-1}(K_1^*)=(\sigma_F^+)^{n-1}(K_0^*)=M$ and $K_0^* \cup K_1^* \subset M$. So $(\sigma_F^+)^{n-1}$ is separated turbulent.

\vphantom{}

Next assume $n\in \mathbb{N}$, $n>1$, and $(\sigma_F^+)^{n-1}:X_F^+ \to X_F^+$ is separated turbulent. Let $K_0$ and $K_1$ be nonempty disjoint closed subsets of $X_F^+$ such that $K_0 \cup K_1 \subset (\sigma_F^+)^{n-1}(K_0) \cap (\sigma_F^+)^{n-1}(K_1)$. Let $m \in \mathbb{N}$ such that for each $l \ge m$, $\pi_{[1,l]}(K_0) \cap \pi_{[1,l]}(K_1)= \emptyset$.

\vphantom{}

Note that $K_0 \cup K_1 \subset (\sigma_F^+)^{n-1}(K_0) \cap (\sigma_F^+)^{n-1}(K_1)$ implies that $K_0 \cup K_1 \subset (\sigma_F^+)^{n-1}(K_0 \cup K_1)  \subset (\sigma_F^+)^{n-1}(K_0) \cup (\sigma_F^+)^{n-1}(K_1) \subset (\sigma_F^+)^{n-1}((\sigma_F^+)^{n-1})(K_0) \cup (\sigma_F^+)^{n-1}(K_1)) \subset (\sigma_F^+)^{2n-2}(K_0) \cup (\sigma_F^+)^{2n-2}$ and, in general, for each $i \ge 1$ $K_0 \cup K_1 \subset (\sigma_F^+)^{in-i}(K_0)\cup (\sigma_F^+)^{in-i}(K_1)$. So, if $in-i \ge m$, then $\pi_{[1,in-i+1]}(K_0)$ and $\pi_{[0,in-i+1]}(K_1)$ are disjoint nonempty closed subsets of $X_F^{in-i}$. Let $A_0=\pi_{[1,in-i+1]}(K_0)$ and $A_1=\pi_{[1,in-i+1]}(K_1)$. Then $A_0$ and $A_1$ are nonempty disjoint closed sets in $X_F^{in-i}$ and $\pi_1(A_0) \cup \pi_1(A_1) \subset \pi_{in-i+1}(A_0) \cap \pi_{in-i+1}(A_1)$. Thus $F$ is CR-turbulent.

\end{proof}

 \vphantom{}
 
 It is known that if $X$ is a compact metric space and $f:X \to X$ is a function on $X$ such that $f^n$ is separated turbulent for some $n$, then $ent(f)>0$. See Block and Coppel \cite{BlockCoppel}, p. 128. It is also easy to see that if $F$ is a CR-turbulent relation on $X$, then $\sigma_F^+:X_F^+ \to X_F^+$ is CR-turbulent, so that $ent(\sigma_F^+)>0$, by noting that, from the proof of the previous theorem, we have $M= \cup_{\alpha \in \Sigma_2} M_{\alpha}$, $(\sigma_F^+)^{n-1}(M)=M$ and $M_{\alpha} \cap M_{\beta} =
 \emptyset $ for $\alpha,\beta \in \Sigma_2$ and $\alpha \ne \beta$. If we define $h:M \to \Sigma_2$ by $h(\mathbf{x})=\alpha$ if and only if $\mathbf{x} \in M_{\alpha}$, then $h \circ (\sigma_F^+)^{n-1}|_M=\sigma_{\Sigma_2} \circ h$. Thus, $\sigma_{\Sigma_2}$ is a factor of $(\sigma^+)^{n-1}|_M$. From this it follows that $ent(\sigma_F^+) >0$. Thus, we have the following result:
 
 \vphantom{}
 
 \begin{theorem} 
 If $X_F^n$ contains two disjoint continua $K$ and $L$ such that $\pi_1(K) \cup \pi_1(L) \subset \pi_n(K) \cap \pi_n(L)$ (i.e., $F$ is CR-turbulent), then $ent(\sigma_F^+)>0$. 

\end{theorem}

\vphantom{}

\begin{exam} 

We know $f:X \to X$ is turbulent implies $f^n:X \to X$ is turbulent (for $n \in \mathbb{N}$). It is natural to wonder whether there is an example of a continuous $f:X \to X$ such that $f^2$ is turbulent but $f$ is not. 

Let $T:[0,1] \to [0,1]$ be the standard tent map defined by 
\[
T(x)= 
\begin{cases}
2x & \text{if } 0 \le x \le 1/2, \\
2-2x & \text{if } 1/2 \le x \le 1. \\
	
\end{cases} 
\]

Define $f:\{0,1, \dots,n-1 \} \times [0,1] \to \{0,1, \dots,n-1 \} \times [0,1]$ by

\[
F(k,x)=
\begin{cases}
(1,T(x)) & \text{if } k=0, \\
((k+1)\mod n,x) & \text{if } k \ne 0. \\
	
\end{cases}
\]
The function $f$ is $n$-turbulent, but not $m$-turbulent for $m < n$. Moreover, if we identify $\{0,1, \dots ,n-1\} \times \{0 \}$, then we get this result on the $n$-od. (See Figures 1 and 2.)

 \end{exam}
 
 \begin{figure}[h!] 
\begin{tikzpicture} 

  \draw [white,very thick] (0,0) -- (0,3);  
  \draw [blue,very thick] (3,0) -- node[very near start, left]{0} (3,3);
  \draw [blue,very thick] (6,0) -- node[very near start, left]{3} (6,3);
  \draw [blue,very thick] (4,0) -- node[very near start, left]{1} (4,3);
  \draw [blue,very thick] (5,0) -- node[very near start, left]{2} (5,3);
  \draw [blue,very thick] (7,0) -- node[very near start, left]{4} (7,3);
  \draw [red,thick] (3.25,1.5) --  (3.75,1.5) -- (3.65,1.6) -- (3.75,1.5) --  (3.65, 1.4);
   
  \node at (3.5,1.8){T};
  \draw [red, thick] (4.25,1.5) -- (4.75,1.5) -- (4.65,1.6) -- (4.75,1.5) -- (4.65,1.4);
  \draw [red, thick] (5.25,1.5) -- (5.75,1.5) -- (5.65,1.6) -- (5.75,1.5) -- (5.65,1.4);
  \node at (4.5,1.8){id};
  \draw [red, thick] (6.25,1.5) -- (6.75,1.5) -- (6.65,1.6) -- (6.75,1.5) -- (6.65,1.4);
  \node at (6.5,1.8){id};
  \node at (5.5,1.8){id};
  
  \draw [red, thick] (7.3,1.4) -- (7.75,1.5);
  
  \draw [red, thick] (7.75,1.5) .. controls (7,6) and (6,3.5) .. (5,3.5);
  \draw [red,thick] (5,3.5) .. controls (4,3.5) and (3,6) ..  (2,3);
  \node at (5,4){id};
  \draw [red, thick] (2,3) .. controls (1,2) and (1.5,1.5) .. (2.5,1.75);
  \draw [red, thick] (2.4,1.85) --(2.5,1.75) -- (2.4,1.65);

  \end{tikzpicture}

\caption{The system $f$ for $n =5$}

\end{figure}
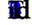

 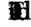
\begin{figure}[h!]
\begin{tikzpicture}
\draw [white,very thick] (0,0) -- (0,3);  
\draw [blue,very thick] (5,3) -- (8,3);
\node at (8.2,3){0};
\draw [blue,very thick] (5,3) -- (5.927,5.853);
\node at (5.99,6.04){1};


\draw [blue,very thick] (5,3) -- (2.573,4.763);
\node at (2.41,4.88){2};

\draw [blue,very thick] (5,3) -- (2.572,1.237);
\node at (2.41,1.12){3};
\draw [blue,very thick] (5,3) -- (5.927,.147);
\node at (5.99,-0.04){4};

\draw[->,red, thick] (6.48,3.26) -- (5.7,4.324);
\node at (6.37,3.99){T};
\draw[->,red, thick] (5.21,4.49) -- (3.96,4.08);
\node at (4.47,4.62){id};
\draw[->,red, thick] (3.65,3.66) -- (3.65,2.34);
\node at (3.3,3){id};
\draw[->,red, thick] (5.21,4.49) -- (3.96,4.08);
\node at (4.47,1.38){id};
\draw[->,red, thick] (3.96,1.92) -- (5.21,1.51);
\node at (6.38,2){id};
\draw[->,red, thick] (5.66,1.65) -- (6.48,2.74);

 \end{tikzpicture}

\caption{The system $f$ when the bottom points are identified}
\end{figure}

\vphantom{}

Next we show that the closed relation $\Gamma(F)$ is CR-turbulent if and only in $f$ is separated $n$-turbulent for some $n$ and then we give a condition under which $f:X \to X$ is turbulent in the classical sense implies $f$ is separated $2$-turbulent. 

\begin{theorem}
Suppose $f:X \to X$ is continuous, and $\Gamma(f)$ is the graph of $f$. Then $\Gamma(f)$ is CR-turbulent if and only if  $f$ is separated $m$-turbulent for some positive integer $m$.  

\end{theorem}

\begin{proof} Suppose $\Gamma(f)$ is CR-turbulent. Then there are $n \in \mathbb{N}$, $n>1$, and disjoint nonempty closed sets $K$ and $L$ in $X_{\Gamma(f)}^{n-1}$ such that $\pi_1(K) \cup \pi_1(L) \subset \pi_n(K) \cap \pi_n(L)$. (Recall that $X_{\Gamma(f)}^{n-1}= \{ \mathbf{x}=(x_1, \dots, x_n):f(x_i)=x_{i+1} \text{ for } 1 \le i <n \}$.)

\vphantom{}

Suppose $\pi_1(K) \cap \pi_1(L) \neq \emptyset$. Then there is some $y \in \pi_1(K) \cap \pi_1(L)$ such that $\mathbf{y} = (y,f(y),f^2(y), \dots, f^{n-1}(y)) \in X_{\Gamma(f)}^{n-1}$. (No other point in $K$ has $1$st coordinate $y$ and no other point in $L$ has $1$st coordinate $y$). Then $\mathbf{y} \in K \cap L$, which is a contradiction. So $\pi_1(K) \cap \pi_1(L) = \emptyset$.

\vphantom{} 

Let $K_0=\pi_1(K)$ and $L_0=\pi_1(L)$. Then $K_0$ and $L_0$ are nonempty disjoint closed subsets of $X$. Furthermore, $K_0 \cup L_0 \subset f^{n-1}(K_0) \cap f^{n-1}(L_0)$: Without loss of generality, suppose $z \in K_0$. Then $(z,f(z), \ldots, f^{n-1}(z)) \in K$ and $f^{n-1}(z) \in \pi_n(K)$. But $\pi_1(K) \cup \pi_1(L) \subset \pi_n(K) \cap \pi_n(L)$, so $z \in \pi_n(K) \cap \pi_n(L) = \{f^{n-1}(x):x \in K \} \cap \{f^{n-1}(x):x \in L \} = f^{n-1}(K_0) \cap f^{n-1}(L_0) $. It follows that $f$ is separated $n-1$-turbulent. Let $m=n-1$. Then $f$ is separated turbulent.

\vphantom{}

Now suppose $f$ is separated $m$-turbulent for some positive integer $m$. Then $f^m$ is separated turbulent and there exist disjoint nonempty closed sets $K$ and $L$ of $X$ such that $K \cup L \subset f^m(K) \cap f^m(L)$. Let $K^+=\{ \mathbf{x}=(x_1, \ldots, x_n) \in X_{\Gamma(f)}^{n-1}:x_1 \in K \}$ and $L^+=\{ \mathbf{x}=(x_1, \ldots, x_n) \in X_{\Gamma(f)}^{n-1}:x_1 \in L \}$. 

\vphantom{}

Suppose $K^+ \cap L^+ \neq \emptyset$. Then there is some point $\mathbf{x}=(x_1, \ldots, x_n) \in K^+ \cap L^+$. Then $x_1 \in K \cap L$, which is a contradiction. Hence, $K^+ \cap L^+ = \emptyset$.

\vphantom{}

Suppose $z \in \pi_1(K^+) \cup \pi_1(L^+)=K \cup L$. Without loss of generality, assume $z \in \pi_1(K^+)$. Then $z \in f^m(K) \cap f^m(L)$, so
there is some $w \in K$ such that $f^m(w)=z$ and there is some $u \in L$ such that $f^m(u)=z$. Then $\mathbf{w}=(w_1, \ldots, w_m) \in X_{\Gamma(f)}^{m-1}$ such that $w_1=w$ and $w_m=z$, and $\mathbf{u}=(u_1, \ldots, u_m) \in X_{\Gamma(f)}^{m-1}$ such that $u_1=u$ and $u_m=z$. Then $z \in \pi_m(K^+) \cap \pi_m(L^+)$ since $\mathbf{w} \in K^+$ and $\mathbf{u} \in L^+$.

\end{proof}


\vphantom{}




The following theorem gives an easy-to-check condition for when a closed relation induces separated turbulence for the second iterate of the shift $\sigma^2$.

\vphantom{}

 \begin{theorem}
 If $X$ is a compact metric space and $F \subset X \times X$ is a closed relation such that $F$ contains closed sets $K$ and $L$ with $\pi_1(K)=\pi_1(L)=X$ and $\pi_1(K \cap L) \cap\pi_2(K \cap L) = \emptyset$, then $\sigma^2_{F^{-1}}$ is separated turbulent.
 \end{theorem}  
 \begin{proof}
 If $(x,y,z) \in (K \star K) \cap (L \star L)$, then $y \in \pi_2(K) \cap \pi_2(L)$ and $y \in \pi_1(K) \cap \pi_1(L)$. So $(x,y) \in K \cap L$ and $(y,z) \in K \cap L$, and it follows that $y \in \pi_2(K \cap L) \cap \pi_1(K \cap L),$ which is a contradiction. Then $K \star K$ and $L \star L$ are disjoint closed sets in $X^3_F$. Now $\pi_1(K \star K)=X$, because if $x \in X$, there is $y$ such that $(x,y) \in K$, which implies that $y \in \pi_2(K) \subset X$, and there is some $z$ such that $(y,z) \in K$. Then $(x,y,z) \in K \star K$ and $x \in \pi_1(K \star K)$. Likewise, $X= \pi_1(L \star L)$. Then, since $\pi_3(K \star K) \subset \pi_1(K \star K)=X$ and $\pi_3(L \star L) \subset \pi_1(L \star L)=X$, $F$ is reverse CR-turbulent for $n=3$, so $F^{-1}$ is CR-turbulent for $n=3$, and the result follows. 
 \end{proof}
 
 \vphantom{}
 
 If $f:I \to I$ is continuous, then turbulence implies separated 2-turbulence. We don't know if this is true for compact metric spaces, but we do have the following theorem whose proof is similar to that of Theorem 3.5.
 
 \vphantom{}
 
 \begin{theorem}
Suppose $X$ is a compact metric space, and $f:X \to X$ is turbulent with $K_0$ and $K_1$ nonempty closed subsets of $X$ such that $K_0 \cap K_1 = \{ x \}$ for some point $x \in X$, and $K_0 \cup K_1 \subset f(K_0) \cap f(K_1)$. Then if $f(x) \neq x$, $f$ is separated 2-turbulent.
\end{theorem}

\begin{proof}
Suppose $f:X \to X$ is a continuous function on a compact metric space $X$, and there are closed sets $K_0,K_1 \subset X$ and $x \in X$ such that $K_0 \cap K_1= \{x \}$ and $K_0 \cup K_1 \subset f(K_0) \cap f(K_1) $ and $f(x) \ne x$. Let $K_0' = K_0 \cap f^{-1}(K_0)$, $K_1' =	K_1 \cap f^{-1}(K_1)$. Note that $K_0'$ and $K_1'$ are closed subsets of $X$. If $K_0' \cap K_1' \ne \emptyset$, then, since $K_0' \subset K_0$ and $K_1' \subset K_1$, $K_0' \cap K_1' =\{x\}$. Thus $x \in f^{-1}(K_0)$ and $x \in f^{-1}(K_1)$, so $f(x) \in K_0 \cap K_1= \{x\}$. But $f(x) \ne x$, so $K_0' \cap K_1' = \emptyset$. If $y \in K_0$, there is $z \in K_0$ such that $f(z)=y$. So $z \in K_0'$ and $f(z)=y$. Thus, $K_0 \subset f(K_0') \subset K_0$ and $f(K_0')=K_0$. Similarly, $f(K_1')=K_1$. It follows that $K_0' \cup K_1' \subset f^2(K_0') \cap f^2(K_1')$. That is, $f$ is $2$-separated turbulent. 
\end{proof}

\vphantom{}

 \vphantom{}
 
 
 \section{Dynamics on finite closed relations}

We use the following theorem from \cite{van}:

\begin{theorem} (This is Theorem 5.6 in \cite{van}.) Let $X$ be a non-empty compact metric space and let $F$ be a finite subset of $X \times X$. The following statements are equivalent.
\begin{enumerate}
	\item $ent(F) \ne 0$.
	\item There are a set $A \subset X$, a positive integer $k$, and an $\epsilon >0$ such that $F$ has a $(k,\epsilon)$-return on $A$. (We do not use this part of the theorem and so do not define the term $(k,\epsilon)$-return.)
	\item There are \begin{enumerate} 
	\item positive integers $k_x$ and $k_y$, 
	\item points $\mathbf{x} \in X_F^{k_x}$ and $\mathbf{y} \in X_F^{k_y}$ such that $\mathbf{x}(k_x) = \mathbf{x}(0)=\mathbf{y}(0)=\mathbf{y}(k_y),$
	\item a positive integer $j$ such that $0 < j \le \min{k_x,k_y}$ and $\mathbf{x}(j)\ne \mathbf{y}(j)$. \end{enumerate}
	\item $X_F^+$ is uncountable.

\end{enumerate}	
\end{theorem}

\vphantom{}

Now assume $ent(F) >0$. Switching to the notation in this paper, part (3) becomes the following: There are (a) positive integers $n_x$ and $n_y$, (b) points $\mathbf{x} \in X_F^{n_x}$ and $\mathbf{y} \in X_F^{n_y}$ such that $\mathbf{x}_{n_x} =\mathbf{x}_{1}=\mathbf{y}_{1} =\textbf{y}_{n_y}$, and (c) a positive integer $1 <j \le \min{n_x,n_y}$ and $\mathbf{x}_{j} \ne \mathbf{y}_{j}$. Let $A=\{\mathbf{x} \star \mathbf{y} \}$ and let $B=\{ \mathbf{y} \star \mathbf{x} \}$. Then $A$ and $B$ are closed subsets of $X_F^{n_x +n_y -1} $ such that $A \cap B = \emptyset$, and 
$$\pi_1(A)=\pi_2(B) = \{\mathbf{x}_1 \}= \{ \mathbf{y}_1 = \pi_{n_x+n_y -1}(A) = \{\pi_{n_x}(\mathbf{x}) \} = \{ \pi_{n_y}(\mathbf{y}) \} = \pi_{n_x + n_y -1}(B).$$ Thus, $\pi_1(A) \cup \pi_1(B) \subseteq \pi_{n_x+n_y -1}(A) \cap \pi_{n_x+n_y-1}(B)$. It follows that $F$ is CR-turbulent. This proves that part(4) of the theorem below implies part(1). 

\vphantom{}
\begin{theorem}
If $x$ is a compact metric space and $F$ is a finite subset of $X \times X$, 
 then the following are equivalent: \begin{enumerate}
\item $F$ is CR-turbulent.
\item There is some $n \in \mathbb{N}$ such that $(\sigma^+_F)^+$ is separated turbulent.
\item There are some $n \in \mathbb{N}$ and two points $\mathbf{x}$ and $\mathbf{y}$ in $X_F^n$ such that $\mathbf{x} \neq \mathbf{y}$ and $\pi_1(\mathbf{x})=\pi_1(\mathbf{y})=
\pi_n(\mathbf{x})=\pi_n(\mathbf{y})$.

\item $ent(\sigma_F^+) \neq 0$.
\item $X_F^+$ is uncountable. 
\item $(X_F,\sigma_F^+)$ is Li-Yorke chaotic.
\item $(X_F,\sigma_F^+)$ has a DC2-scrambled Cantor set.

\end{enumerate}
\end{theorem}

\begin{proof} That (1) is equivalent to (2) is proved Theorem 3.1. In the theorem from \cite{van} above,  it is shown that (3) - (7) are equivalent (Theorem 5.7 in \cite{van} and Corollary 7 in \cite{van}.) Clearly, (3) implies (1). Hence, all the statements are equivalent.

\end{proof}

\vphantom{} 


\section{Closed relations that are neither CR-turbulent nor reverse CR-turbulent but do have positive entropy}


Results of M. Misiurewicz, in \cite{Misiurewicz}, imply that a continuous map $f$ on a closed interval has positive topological entropy if and only if $f^n$  separated turbulent for some $n$. This gives rise to the following question.

\vphantom{}

\textbf{Question.} If $F \subset [0,1]^2$ is a closed relation that is the union of the graphs of two continuous functions defined on $[0,1]$, does $ent(F)>0$ imply that $F$ is CR-turbulent or reverse CR-turbulent?

\vphantom{}

Below we give some results to this question and related examples.

\vphantom{}

\begin{lemma} Suppose $X$ is a compact metric space and $G$ is a closed subset of $X \times X$ and $a,b \in X$ such that $a \ne b$ and 
\begin{enumerate}
\item if $(x,a) \in G$, then $x=a$;
\item if $(b,x) \in G$, then $x=b$; and 
\item if $(s,u)\in G$ and $(t,u) \in G$ and $s \ne t$,	then $\{s,t,u\} \cap \{a,b\} \ne \emptyset$.
\end{enumerate} 
 Then $G$ is not CR-turbulent.

\end{lemma}

\begin{proof} Suppose $A$ and $B$ are nonempty disjoint closed sets in $X_G^+$ and $n \in \mathbb{N}$ such that $A \cup B \subset (\sigma^+_G)^n(A) \cap (\sigma^+_G)^n(B)$. Let $\mathbf{z} \in A$. Then there is $\mathbf{x} \in A$ and there is $\mathbf{y} \in B$ such that $(\sigma^+_G)^n(\mathbf{x})=(\sigma^+_G)^n(\mathbf{y})=\mathbf{z}$. But if there is $i \in \mathbb{N}$ such that $z_i=a$, then $x_j=y_j=a$ for each $j \le i$. So $\mathbf{x}=\mathbf{y}$. It follows that for each $\mathbf{w}=(w_1,w_2, \ldots) \in A$ and for each $i \in \mathbb{N}$, $w_i \ne a$. Similarly, for each $\mathbf{w} \in B$ and for each $i \in \mathbb{N}$, $w_i \ne a$. Since $\mathbf{x} \ne \mathbf{y}$ and $(\sigma^+_G)^n(\mathbf{x})=(\sigma^+_G)^n(\mathbf{y})$, there is $i \le n$ such that $x_i \ne y_i$ and $x_{i+1}=y_{i+1}$. By (3), $\{x_i,y_i,x_{i+1}\} \cap \{a,b \} \ne \emptyset$. Thus, $b \in \{x_i,y_i,x_{i+1}\}$, so by (2), $\mathbf{z}=\overline{\mathbf{b}}=(b,b,b,b,\dots)$. Thus, $A = \{ \overline{\mathbf{b}} \}$ and similarly, $B = \{ \overline{\mathbf{b}} \}$, which contradicts $A \cap B = \emptyset$. It follows that $G$ is not CR-turbulent.
	
\end{proof}

\begin{theorem}
Suppose $X$ is a compact metric space and $f:K_1 \to K_2$ is a homeomorphism from closed set $K_1 \subset X$ to closed set $K_2 \subset X$. Suppose there are $a,b \in X$ such that $a \ne b$ and $(\{(a,x):x \in X \} \cup \{(x,b):b \in X\}) \cap \{(x,f(x)):x \in K_1 \} \subset \{(x,x):x \in X \}$. If $G=\{(a,x):x \in X \} \cup \{(x,b):b \in X \} \cup \{(x,f(x)):x \in K_1\}$, then $\pi_1(G)=\pi_2(G) = X$; $G$ is not CR-turbulent; and $G$ is not reverse CR-turbulent. 	
\end{theorem}

\begin{proof} Assume $X,f,a,b$ are as above. It is easy to verify the following:

\begin{enumerate}
	\item If $(x,a) \in G$, then $x=a$.
	\item If $(b,x) \in G$, then $x=b$.
	\item If $(s,u) \in G$ and $(t,u) \in G$ and $s \ne t$, then $\{s,t,u \} \cap \{a,b\} \ne \emptyset$.
\end{enumerate}
So, by Lemma 5.1, $G$ is not CR-turbulent. The proof that $G^{-1}$ is not CR-turbulent is the same but the roles of $a$ and $b$ are switched.
	
\end{proof}

\begin{exam} A closed subset $F$ of $[0,1]^2$ such that $\pi_1(F)=\pi_2(F)=[0,1]$, $ent(F)>0$, and $F$ is neither CR-turbulent nor reverse CR-turbulent: Let $f$ be a homeomorphism from the embedding of the Cantor middle-third set in $[0,1]$ onto the same embedding of the middle -third Cantor set in $[0,1]$ with $ent(f) >0$. Choose $a,b \in [0,1]$ such that $a \ne b$ and neither $a$ nor $b$ is in the Cantor middle-third set in $[0,1]$. Let $V_a=\{(a,y)|y \in [0,1] \}$ and $H_b= \{(x,b)|x \in [0,1] \}$. If $F=\Gamma(f) \cup V_a \cup H_b$, then $ent(F)>0$, $\pi_1(F)=\pi_2(F)= [0,1]$, and $F$ is not CR-turbulent and not reverse CR-turbulent. (See Figure 3.)
	
\end{exam}

\begin{figure}[h!]
\begin{tikzpicture}

\draw [white,very thick] (0,0) -- (0,3);  

\draw (1,0) -- (7,0) -- (7,6) -- (1,6) -- (1,0);

\draw [blue,very thick] (2.7,0) -- (2.7,6);
\draw [blue,very thick] (1,5) -- (7,5);

\filldraw[blue] (1.5,1) circle (2pt); 
\filldraw[blue] (2.2,3) circle (2pt); 
\filldraw[blue] (1.5,0) circle (2pt);
\filldraw[blue] (1,1) circle (2pt);
\filldraw[blue] (2,0) circle (2pt);\filldraw[blue] (2.2,0) circle (2pt);  
\filldraw[blue] (1,3) circle (2pt); 
\filldraw[blue] (2.2,5.3) circle (2pt);   
\filldraw[blue] (1,5.3) circle (2pt);
\filldraw[blue] (1,1.2) circle (2pt);  
\filldraw[blue] (6.3,0) circle (2pt); 
\filldraw[blue] (5,1.2) circle (2pt); 
\filldraw[blue] (5,0) circle (2pt);
\filldraw[blue] (3,1.8) circle (2pt); 
\filldraw[blue] (3,0) circle (2pt);  \filldraw[blue] (1,1.9) circle (2pt); 
\filldraw[blue] (2.9,5.8) circle (2pt);
\filldraw[blue] (2.9,0) circle (2pt); 
\filldraw[blue] (1,5.8) circle (2pt); 
\filldraw[blue] (6.8,4.6) circle (2pt); 
\filldraw[blue] (6.8,0) circle (2pt); 
\filldraw[blue] (1,4.6) circle (2pt); 
\filldraw[blue] (5.6,2.5) circle (2pt); 
\filldraw[blue] (5.6,0) circle (2pt);  
\filldraw[blue] (1,2.5) circle (2pt); 
\filldraw[blue] (2.5,4.1) circle (2pt); 
\filldraw[blue] (2.5,0) circle (2pt); 
\filldraw[blue] (1,4.1) circle (2pt); 
\filldraw[blue] (4.1,4.3) circle (2pt); 
\filldraw[blue] (4.1,0) circle (2pt); 
\filldraw[blue] (1,4.3) circle (2pt);

\end{tikzpicture}
\caption{Illustration of a closed relation with satisfying the conditions of Example 5.3 which has no CR-turbulence and no reverse CR-turbulence. Blue dots on the axes represent a Cantor set in the interval. Blue dots in the interior represent a the graph of a homeomorphism from that Cantor set to itself.}
\end{figure}

\begin{exam}
Let $CF$ denote the Cantor Fan. There is a homeomorphism of the Cantor Fan CF with positive entropy. (The endpoints of $CF$ form a Cantor set. Then there is a homeomorphism $h$ with positive entropy on the endpoints $E(CF)$. Extend this homeomorphism to a homeomorphism $h'$ on $CF$ by defining that a leg $A$ in $CF$ is moved to the leg $B$ if $h(e(A))=e(B)$, and defining $h'(x)$ in the obvious way.) Then $\Gamma(F) \subset CF \times CF$ is the graph of such a homeomorphism, and $\Gamma(F)$ is neither CR turbulent nor reverse CR turbulent.  
\end{exam}

\vphantom{}

For closed relations $F$ on compact metric spaces $X$, we hope to use the results in this section in the future to explore further the connection between $ent(F)>0$ and the CR-turbulence of $F$.  

 \section{Turbulence on a large but simple class of closed relations on $[0,1]$}

The result discussed below gives a large class of closed relations on an interval that are CR-turbulent. 

\vphantom{}

\begin{figure}[h] 
\begin{tikzpicture}

\begin{axis}[axis lines=middle, xlabel={$x$}, ylabel={$y$}, domain=0:1, xmin=0,xmax=1, ymin=0, ymax=1, samples=100, legend pos=north west]

\addplot[blue,thick]{2*sqrt(x)};
\addplot[red, thick]{.5 *x^2};
\end{axis}

\end{tikzpicture}

\caption{The maps $f$ and $g$}

\end{figure}

The remaining results in this section assume the following:

 Suppose $0<a<1$ and $f:[0,a] \to [0,1]$ is a continuous injection such that 
\begin{itemize}
\item $f(0)=0$,
\item $f(a)=1$, and 
\item $f(x)>x$ for each $a \ge x >0$.
\end{itemize}

Suppose $0<b<1$ and $g:[0,1] \to [0,1]$ is a continuous injection such that 
\begin{itemize}
\item $g(0)=0$,
\item $g(1)=b$, and 
\item $g(x)<x$ for each $1 \ge x >0$.
\end{itemize}

Note that \begin{enumerate}
\item $f^{-n}(1) < f^{-n+1}(1)$ for each $n \in \mathbb{N}$,
\item $\lim_{n \to \infty} f^{-n}(1)=0$,
\item $g^n(1)<g^{n-1}(1)$ for each $n \in \mathbb{N}$, and
\item $\lim_{n \to \infty}g(1)=0$.
\end{enumerate}







\begin{prop}

Let $\alpha$ be such that $\alpha =f^{-1}\circ g(1)$. Thus $f(\alpha)=g(1)$. Note that $\alpha <a$. Let $f \cup g$ be the closed relation defined by $(f \cup g)(x)= \{f(x),g(x)\}.$ Then $$(f \cup g) ([\alpha,1])=g([\alpha,1]) \cup f([\alpha,a])$$ $$=[g(\alpha),g(1)] \cup [f(\alpha),f(a)]$$ $$=[g(\alpha),g(1)] \cup [g(1),1]$$ $$=[g(\alpha),1]. $$ 
Furthermore, this holds for $0< \beta \le \alpha$, i.e., 
$$(f \cup g) ([\beta,1])=g([\beta,1]) \cup f([\beta,a])$$ $$=[g(\beta),g(1)] \cup [f(\beta),f(a)]$$ $$=[g(\beta),g(1)] \cup [g(1),1]$$ $$=[g(\beta),1], $$
and for each $n$, $$(f \cup g)^n([\alpha, 1])=[g^n(\alpha),1],$$
and  $$(f \cup g)^n([\beta, 1])=[g^n(\beta),1].$$
\end{prop}

\begin{prop} Given $p \in (0,1)$, there are $z=z(p) \in (0,p)$ and $M=M(p)$ such that $(f \cup g)^M([z,p])=[g^M(z),1] \supset [\alpha,1]$. 

\end{prop}

 \begin{proof}
Note that $f^{-1}\circ g(p) \le \alpha$.

\textbf{Case 1: $a \le p$.} Let $z$ be such that $f(z)=g(p)$. Then $z=f^{-1}\circ g(p) \le \alpha$. Let $M=1$. Then $(f \cup g)([z,p]) =f([z,a]) \cup g([z,p]) =[g(p),1] \cup [g(z),g(p)] = [g(z),1] \supseteq [\alpha,1]$. 

\textbf{Case 2: $p <a$.} Let $m \in \mathbb{Z}$ such that $f^{m-2}(p) <a \le f^{m-1}(p)$. Let $z=f^{-m} \circ g^m(p) \le \alpha$. Then $f^{m-1}([z,p]) \supseteq [f^{m-1}(z),a]$, so $f^m([z,p]) \supseteq f([f^{m-1}(z),a])= [f^m(z),f(a)]=[g^m(p),1]$. Now $g^m([z,p])= [g^m(z),g^m(p)] = [g^m(z),f^m(z)]$. Let $M=m$. Then \newline $(f \cup g)^M([z,p]) \supseteq [g^M(z),1] \supseteq [\alpha,1]$. 

\end{proof}

\vphantom{}

\begin{theorem} The closed relation $f \cup g$ defined above is CR-turbulent and reverse CR-turbulent. 
\end{theorem}

\begin{proof}
Let $\alpha= f^{-1}\circ g(1)$ (so $f(\alpha)=g(1)=b$). Let $p \in (0,1)$ and let $M=M(p)$ and $z=z(p)$ be defined as in Proposition 6.2. Let $\beta=g^{-1}(z)$, $\gamma=f^{-1}(\alpha)$. Then notice that $f([\gamma,a])=[\alpha,1]$. Since $g^n(z) \to 0$ and $g^n(\alpha) \to 0$ as $n \to \infty$, there is $N \in \mathbb{N}_0$ such that $g^{M+N}(z) < \min \{ \gamma, \beta \}$ and $g^{M+N}(\alpha) < \min \{\gamma, \beta \}$.

\vphantom{}

Let $A$ and $B$ be subsets of $I_{f \cup g}^{N+M+1}$ defined by 

$$A= \{ (x_1,x_2, \ldots, x_{N+M+1}): x_1 \in [\beta,1] ,g(x_1)=x_2, x_3 \in (f \cup g)  (x_2), \ldots, x_{N+M+1} \in(f \cup g)(x_{N+M}) \},$$ and $$B= \{ (x_1,x_2, \ldots, x_{N+M+1}):x_1 \in [\gamma,\alpha] ,f(x_1)=x_2, x_3 \in (f \cup g)  (x_2), \ldots, x_{N+M+1} \in(f \cup g)(x_{N+M}) \}.$$

Now, if $(x_1,x_2,\ldots,x_{N+M+1}) \in A \cap B$, then $x_1 \in [\beta,1] \cap [\gamma,\alpha]$, and $x_2=g(x_1)=f(x_1)$. But this cannot be since $g(x_1) < f(x_1)$. Then $A \cap B = \emptyset$. Let $$A'= \{ (x_1,x_2, \ldots, x_{N+M+1}): x_1 \in [\beta,1] ,g(x_1)=x_2, x_3 \in g(x_2), \ldots, x_{N+M+1} \in g(x_{N+M}) \},$$ and let $$B'= \{ (x_1,x_2, \ldots, x_{N+M+1}):x_1 \in [\gamma,\alpha] ,f(x_1)=x_2, x_3 \in  g  (x_2), \ldots, x_{N+M+1} \in  g(x_{N+M}) \}.$$ Then $\pi_1(A) =[\beta,1]$,  because $[\beta,1] = \pi_1(A') \subset \pi_1(A)$. Likewise, $[\gamma,\alpha]=\pi_1(B') \subset \pi_1(B)$, so $\pi_1(B)=[\gamma.\alpha]$.

\vphantom{}

Since $g^{M+N}(z) < \beta$ and $z \le \alpha$, it follows from Proposition 6.1 that $$(f \cup g)^n([z, 1])=[g^n(z),1].$$ Since $g^{M+N}(z) < \gamma$, it follows from Proposition 6.2, that $[g^{M+N}(z),1] \supset [g^M(z),1] \supset [\alpha,1]$, so $\pi_{M+N+1}(A)=[g^{M+N}(z),1] \supseteq [\beta,1] \cup [\gamma,a]$. A similar argument yields the result that $\pi_{M+N+1}(B)=[g^{M+N}(\alpha),1] \supseteq [\beta,1] \cup [\gamma, a]$. 

\vphantom{}

Hence, $f \cup g$ is CR-turbulent. By considering the same argument applied to $f^{-1}$ and $g^{-1}$, it follows that $f \cup g$ is also reverse CR-turbulent. 

\end{proof}


\vphantom{}

  A continuous turbulent function on a closed interval must have periodic points of all periods (see \cite{BlockCoppel}, Chapter II, Lemma 3). However, the following example illustrates a striking contrast. In this example of a turbulent continuous function on a Lelek fan, the associated shift map (which is continuous, surjective, and turbulent) has exactly one periodic point: the fixed point (0, 0, 0, \ldots).

\vphantom{}

\begin{exam}
A turbulent map on the Lelek fan with one fixed point and no other periodic orbits:  Let $a=1/3$, $b=2$, and let $F=\{(x,ax):x \in [0,1] \} \cup \{ (x,bx): 0 \le x \le 1/2 \}$. Then $\sigma_L: I_F^+ \to I_F^+$ is both turbulent and CR-turbulent, but has no periodic points other than $(0,0,\dots)$, which is a fixed point. That  this example is a Lelek fan is proved in  \cite{banic1}. That it has no periodic points other than the one fixed point is proved in \cite{BE2}.
\end{exam}

\vphantom{}

\textbf{Observation.} It follows from the theorem above that $F_{ab}$ is CR-turbulent and reverse CR-turbulent for all $a,b$ such that $0 <a<1<b$.

\vphantom{}

\textbf{Observation.} If $F \subset G \subset X \times X$ and $F$ is CR-turbulent, then $G$ is CR-turbulent.

\vphantom{}

\textbf{Question.} If $f$ is a closed subset of $[0,1]^2$ such that $F$ is the union of the graphs of two continuous functions on $[0,1]$ and $ent(F)>0$, must $F$ be CR-turbulent or reverse CR-turbulent?

\section{Separated, continuum-wise semi-turbulent homeomorphisms on smooth fans}

\vphantom{}

While homeomorphisms cannot be turbulent, reverse turbulent, or separated turbulent, they can satisfy a related condition called ``separated continuum-wise semi-turbulence'':

\vphantom{}

\begin{definition}

Let $X$ be a compact metric space and $h:X \to X$ be a continuous function. We say that $h$ is \emph{semi-turbulent} if there exist a continuum $Y$, disjoint closed sets $A,B$ in $X$ and continuous functions $\phi:X \to Y$, $f:Y \to Y$, such that 

\begin{enumerate}
	\item $f \circ \phi = \phi \circ h$,
	\item$\phi(A) \cap \phi(B) = \emptyset$, and 

	\item $\phi(A) \cup \phi(B) \subset \phi(h(A)) \cap \phi(h(B))$.

\end{enumerate}

See Figure 5.

\end{definition}

\begin{definition}
Let $X$ be a compact metric space and $h:X \to X$ be a continuous function. We say that $h$ is \emph{separated, continuum-wise semi-turbulent} if there exist a continuum $Y$, disjoint subcontinua $A,B$ in $X$ and continuous functions $\phi:X \to Y$, $f:Y \to Y$, such that 

\begin{enumerate}
\item $f \circ \phi = \phi \circ h$,

\item $\phi(A) \cap \phi(B) = \emptyset$, and 

\item $\phi(A) \cup \phi(B) \subset \phi(h(A)) \cap \phi(h(B))$.
\end{enumerate}

Note that this implies that 
$$f \circ \phi(A) \cap f \circ \phi(B) \supset \phi(A) \cup \phi(B),$$
and hence $f$ is separated, continuum-wise turbulent.

\end{definition}

\vphantom{}

\begin{figure}[h!]
\centering
\noindent \begin{tikzpicture}[node distance=1.5cm, auto]
  \node (X1) {$X$};
  \node (X2) [right of=X1] {};
   \node (X3) [right of=X2] {$X$};
    \node (Z1) [below of=X1] {};
  \draw[->] (X1) to node {$h$} (X3);
  \node (Y1) [below of=Z1] {$Y$};
  \node (Y2) [right of=Y1] {};
    \node (Y3) [right of=Y2] {$Y$};
   \draw[->] (Y1) to node {$f$} (Y3);
      \draw[<-] (Y1) to node {$\phi$} (X1);
            \draw[<-] (Y3) to node {$\phi$} (X3);
     \end{tikzpicture}
\caption{The diagram} \label{diagram1}
\end{figure}

\begin{figure}

\begin{tikzpicture}

\draw (6,0) -- (9,0) -- (9,3) -- (6,3) -- (6,0);

\draw [blue, very thick] (6,0) -- (7,3) -- (8,0) -- (9,3);
	
\end{tikzpicture}

\caption{3-fold map on [0,1]}

\end{figure}

\textbf{Note:} Continuum-wise semi-turbulence implies semi-conjugacy to turbulence. Thus, for the homeomorphism $h$ above, $ent(h) >0$.

\begin{exam} Let $F$ be the 3-fold map on the interval $I=[0,1]$ defined by $F(x)=3x \text{ if } 0 \le x \le 1/3$, $F(x)= -3x+2$ if $1/3 \le x \le 2/3$, $F(x)= 3x-2$ if $2/3 \le x \le 1$, and let $X= \underleftarrow{\lim} (I,F)$. Let $\sigma:X \to X$ be the shift homeomorphism, and let $\pi_1:X \to I$ be the projection of $X$ to its first component. Let $A$ be any component of $\pi^{-1}_1([0,1/3])$ and let $B$ be any component of $\pi^{-1}_1([2/3,1])$. Then $\pi_1(A)=[0,1/3]$, $\pi_1(B)=[2/3,1]$, $\pi_1 \circ \sigma = F \circ \pi_1$, and $\pi_1 \circ \sigma (A) = [0,1] = F \circ \pi_1(B)$. Hence, $\sigma$ is separated, continuum-wise, semi-turbulent. However, since $\sigma$ is a homeomorphism, $\sigma$ is not separated, continuum-wise turbulent. See Figure 6.

\end{exam}

\begin{definition} The following is an equivalent definition of the Cantor fan given in Section 2. It is the definition we need here:
	Let $C$ be a Cantor set and $F_p$ be the quotient of $C\times [0,1]$ formed by identifying $\{(c,0)\}_{c\in C}$ with a point $p$. Any continuum homeomorphic to $F_p$ is called a \emph{Cantor fan} and the point $p$ is called the \emph{top} of the Cantor fan.
\end{definition}

\begin{definition}
	Recall that a continuum $Y$ that can be embedded in $F_p$ is called a \emph{ smooth fan.} We do not discuss non-smooth fans, so we define them as follows: A \emph{non-smooth fan} is one that cannot be embedded in the Cantor fan. (It is possible to give a more direct definition, but it is rather involved, and we do not need it here.)
\end{definition}

{\bf Note:} Arcs are trivial smooth fans and any non-trivial smooth fan that is embedded in $F_p$ must contain the top $p$.\\

\vphantom{}

\textbf{Notation.} If $X \subset Y \times Z$, then $\pi_1: X \to Y$ and we define $\pi_1(y,z)=y$ for $(y,z) \in X$. Thus, for $y \in Y$, $\pi_1^{-1}(y)=\{ (y,z):z \in X \}$. 
 
\vphantom{}

\begin{definition}

Let $C$ be a Cantor set and $S$ be a closed subset of $C$. Then $B^S$ is a \textit{blade space} if $B^S$ is a compact subset of $S \times [0,1]$ such that $\pi^{-1}_1(c)$ is connected and contains $(c,0)$ for each $c \in S$. Also, for each $c \in S$, there is $\rho(c) \in [0,1]$ such that $B^S \cap (\{c \} \times [0,1])= \{(c,t): t \le \rho(c) \}$. (Note: it is possible that $\rho(c)=0$.)

\end{definition} 

\vphantom{}

Notice that if $B^S\subset S\times [0,1]$ is a blade space, then the identification of $\{(c,0)\}_{c\in S}$ with a point $p$, $B^S_p$ is a smooth fan. Let $\nu_S:B^S\longrightarrow B^S_p$ be the natural map.

\begin{prop}
	If $Y$ is a smooth fan, then there is a blade space $B^S$ such that $Y$ is homeomorphic to $B^S_p$. 
\end{prop}
\begin{proof}
If $Y$ is a arc then $Y$ is homeomorphic to $\{c\}\times [0,1]$ for any $c\in C$, which is a blade space. So assume that $Y$ is a non-trivial smooth fan.
Let $\nu_C:C\times [0,1]\longrightarrow F_p$ be the natural map. Then $\nu_C^{-1}(Y)$ is a blade space since $C\times\{0\}=\nu_C^{-1}(p)$. 
\end{proof}

\begin{definition} Let $Y$ be a continuum and $A,B$ be disjoint closed subsets of $Y$. Let $I$ be an arc and $\varphi:I \to Y$ be a map. An $(A,B)$\emph{-zigzag set for} $\varphi$ is a subset $\{x_i\}_{i=1}^n \subset I$ such that $x_1 < x_2 < \cdots <x_n$ and $\varphi(x_i) \in A$ if $i$ is even, and $\varphi(x_i) \in B$ if $i$ is odd. (Or, similarly, $\varphi(x_i) \in B$ if $i$ is even, and $\varphi(x_i) \in A$ if $i$ is odd.) The $(A,B)$\emph{-zigzag number for} $\varphi$ is the maximum cardinality that an $(A,B)$-zigzag set can have. (A proof that this maximum cardinality exists is given below.)

\end{definition}

\vphantom{}

\begin{prop} There is $M \in \mathbb{N}$ such that if $\{x_i \}_{i=1}^n$ is an $(A,B)$-zigzag set for $\varphi$, then $n \le M$.
\end{prop}

\begin{proof} Let $\epsilon= \frac{d(A,B)}{2}$. Since $\varphi$ is uniformly continuous, there is $\delta >0$ such that if $x,y \in I$ and $d(x,y)< \delta$, then $d(\varphi(x),\varphi(y)) < \epsilon$. Let $M = \ulcorner \frac{diam (I)}{\delta} \urcorner$. Thus, $M$ is the ``ceiling function'' or the ``round up to the nearest positive integer function''. Then it follows from the pigeon-hole principle that if $\{x_i \}_{i=1}^n$ is an $(A,B)$-zigzag set for $\varphi$, then $n \le M$.

\end{proof}

\begin{lemma} Let $X$ be a blade space, $Y$ be a continuum, $A,B$ be disjoint closed subsets of $Y$, and $\phi:X \to Y$ be continuous (where $X \subset C \times [0,1]$).  For each $c \in C$, let $M(c)$ be the $(A,B)$-zigzag number for $\phi|_{\{c\} \times [0,\rho(c)]}: \{c \} \times [0,\rho(c)] \to Y$. Then $\{ M(c):c \in C \}$ is bounded. 
\end{lemma}

\begin{proof}
Since $X$ is compact, $\phi$ is uniformly continuous. Let $\epsilon = \frac{d(A,B)}{2}$, and let $\delta>0$ such that if $x,y \in X$ and $d(x,y) < \delta$, then $d(\phi(x), \phi(y))< \epsilon$. Since $\pi^{-1}_1(c) \subset \{c \} \times [0,1]$, it follows that $M(c) \le \ulcorner \frac{1}{\delta} \urcorner$, for each $c \in C$.
\end{proof}

\vphantom{}

\textbf{Notation.} In the case that the lemma above holds, let $M_{\phi(A,B)} := \max \{M(c):c \in C \}$ denote the $(A,B)$-zigzag number for $\phi$.

\vphantom{}

Let $\{a_i \}_{i=1}^{2^n} \subset I$ such that $a_i <a_{i+1}$ for each $i$. Let $C_n= \{ (\alpha_1, \dots, \alpha_n): \alpha_i \in \{0,1 \} \}$. Let $\Psi:\{a_i \}_{i=1}^{2^n} \to C_n$ be one-to-one and onto. For each $k \in \{1, \dots,n \}$, let $\{i_j^k\}_{j=1}^{n_k}$ ($n_k \le 2^n$) be a strictly increasing sequence in $\mathbb{N}$ of maximal cardinality such that 
$$ \pi_k \circ \Psi(a_{i_j^k}) \ne \pi_k \circ \Psi(a_{i_j^k+1}).$$ (\textit{Maximal} means that if $i \notin \{i_j^k\}_{j=1}^{n_k}$, then $\pi_k \circ \Psi(a_{i}) = \pi_k \circ \Psi(a_{i+1})$.)  

\begin{lemma}
Given the above, there exists a $k \in \{1, \ldots,n \}$ such that $n_k \ge \frac{2^n-1}{n}$.

\end{lemma} 

\begin{proof}
For each $i \in \{1,2,...,2^n-1\}$ there is $k_i \in \{1,\dots,n\}$ such that $\pi_{k_i} \circ \Psi(a_i) \ne \pi_{k_i} \circ \Psi(a_{i+1})$, since $ \Psi(a_i) \ne \Psi(a_{i+1})$. If $n_k < \frac{2^n-1}{n}$ for each $k \in \{1, \ldots,n \}$, then it follows that there are less than $n(\frac{2^n-1}{n})=2^n-1$ distinct $i$'s such that $ \Psi(a_i) \ne \Psi(a_{i+1})$, so it follows from the pigeon-hole principle that there is at least one $i \in \{1,2,...,2^n-1\}$ such that $ \Psi(a_i) = \Psi(a_{i+1})$, which is a contradiction.

\end{proof}

\vphantom{}

\begin{definition}

 Let $Z \subset X \times Y$ such that $\pi_1(Z)=X$. Then the function $g:Z \to Z$ is a \emph{fiber-injection} if there exist a continuous function $\mu:X \to X$ and, for each $x \in X$, a continuous one-to-one function $h_x: \pi_1^{-1}(x) \to \pi_1^{-1}(\mu(x))$ such that $g(x,y) = (\mu(x), h_x(y))$ and $g$ is continuous.

\end{definition} 

\vphantom{}

\begin{lemma} Let $X$ be a blade space and $g:X \to X$, where $g(c,x)=(\mu(c),h_c(x))$ and $h_c(0)=0$ for each $c \in C$, be a fiber-injection. Then $g$ is not separated, continuum-wise semi-turbulent.
\end{lemma}

\vphantom{} 


\begin{proof} Suppose to the contrary that there are disjoint subcontinua $A,B$ in $X$, a continuum $Y$ and maps $\phi:X \to Y$, $f:Y \to Y$ such that 
\begin{enumerate}
\item $f \circ \phi = \phi \circ g$,
\item $\phi(A) \cap \phi(B) = \emptyset$, and
\item $\phi(A) \cup \phi(B) \subset \phi \circ g(A) \cap \phi \circ g(B)$.

\end{enumerate}

\vphantom{}

Let $M$ be the $(A,B)$-zigzag number for $\phi$. Let $c_A \in C$ such that $A \subset \pi_1^{-1}(c_A)$. Then it follows that 
$$g^n(A) \subset \pi_1^{-1}(\mu^n(c_A)) \subset \{\mu^n(c_A) \} \times [0,1]$$ for each $n$.

\vphantom{}

Since $\phi \circ g(A) \supset \phi(A) \cup \phi(B)$, let 
$$A_{(0,0)}=g(A) \cap \phi^{-1}(\phi(A))$$ and $$A_{(0,1)}=g(A) \cap \phi^{-1}(\phi(B)).$$

Now $\phi(A_{(0,0)})=\phi(g(A) \cap \phi^{-1}(\phi(A))) \subset \phi(g(A)) \cap \phi \circ \phi^{-1} \circ \phi(A)=\phi(g(A)) \cap \phi(A) = \phi(A)$, so $\phi(A_{(0,0)}) \subset \phi(A)$. Also, $\phi(A) \subset \phi \circ g(A)$, and if $y \in \phi(A)$, then $y \in \phi \circ g(A)$, so there is $x \in g(A)$ such that $\phi(x)=y \in \phi(A)$, and it follows that $x \in \phi^{-1}(y) \subset \phi^{-1}(\phi(A))$. Hence, $x \in g(A) \cap \phi^{-1}(\phi(A))=A_{(0,0)}$, so $\phi(A) \subset \phi(A_{(0,0)})$. Then  
 $\phi(A_{(0,0)})=\phi(A)$. The proof that 
$\phi(A_{(0,1)})=\phi(B)$ is similar, and we omit it. Furthermore, $A_{(0,0)} \cap A_{(0,1)} = \emptyset$ (since $\phi(A) \cap \phi(B) = \emptyset$).

\vphantom{}

Next, since $f \circ \phi(A) \cap f \circ \phi(B) \supset \phi(A) \cup \phi(B)$, it follows that  
$$ \phi \circ g(A_{(0,0)})=f \circ \phi(A_{(0,0)})=f \circ \phi(A) \supset \phi(A) \cup \phi(B)$$   and $$ \phi \circ g(A_{(0,1)})=f \circ \phi(A_{(0,1)})=f \circ \phi(B) \supset \phi(A) \cup \phi(B).$$ Continuing inductively, let 

$$A_{(p_1,p_2, \dots,p_n,0)}=g(A_{(p_1, \ldots,p_n)}) \cap \phi^{-1}(\phi(A))$$ and $$A_{(p_1,p_2, \dots,p_n,1)}=g(A_{(p_1, \ldots,p_n)}) \cap \phi^{-1}(\phi(B))$$ where $p_i \in \{0,1\}$. Note that 
$$\phi \circ g(A_{(p_1,p_2, \dots,p_n,0)}) \supset \phi(A) \cup \phi(B),$$ $$\phi \circ g(A_{(p_1,p_2, \dots,p_n,1)}) \supset \phi(A) \cup \phi(B),$$ and
$$A_{(p_1,p_2, \dots,p_n,0)}  \cap A_{(p_1,p_2, \dots,p_n,1)}= \emptyset.$$

\vphantom{}

Let $M$ be the $(\phi(A),\phi(B))$-zigzag number for $\phi$ and choose $n$ such that $M<\frac{2^n-1}{n}$. Notice that $A_{(p_1,p_2, \dots,p_n)} \subset \{ \mu^{n-1}(c_A)\} \times \{0,1 \}$ for each $(p_1,p_2, \dots,p_n) \in C_n$. Since $\{A_{(p_1,p_2, \dots,p_n)}\}_{(p_1,p_2, \dots,p_n) \in C_n}$ are all disjoint, we can pick $\{x_i \}_{i=1}^{2^n} \subset \{\mu^{n-1}(c_A)\} \times \{0,1 \}$ such that $x_i <x_{i+1}$ and each $x_i$ is in a distinct element of  $\{A_{(p_1,p_2, \dots,p_n)}\}_{(p_1,p_2, \dots,p_n) \in C_n}$. Let $\Psi:\{x_i \}_{i=1}^{2^n} \to C_n$ be defined by $\Psi(x_i)=(p_1,p_2, \dots,p_n)$ if and only if $x_i \in A_{(p_1,p_2, \dots,p_n)}$. Then by the previous lemma, there is $k$ and a strictly increasing sequence of natural numbers $\{i_j^k\}_{i=1}^{n_k}$ such that $n_k \ge \frac{2^n-1}{n}>M$ such that $\pi_k \circ \Psi(x_{i_j^k}) \neq \pi_k \circ \Psi(x_j^{k+1})$. Then since $g$ is one-to-one and order preserving, $g^{k-n}(x_{i_j^k}) < g^{k-n}(x_{i_{j+1}^k})$, $\phi \circ g^{k-n}(x_{i_j^k}) \in \phi(A) \cup \phi(B)$, and 
$\phi \circ g^{k-n}(x_{i_j^k}) \in \phi(A)$ if and only if $\phi \circ g^{k-n}(x_{i_{j+1}^k}) \in \phi(B)$. Hence, $\{g^{k-n}(x_{i_j^k})\}_{j=1}^{n_k}$ is a $(\phi(A),\phi(B))$-zigzag set for $\phi$ with cardinality $n_k>M$, which is a contradiction. 
\end{proof}

\begin{definition}
Let $Z\subset X\times Y$ be such that $\pi_1(Z)=X$ and let $y\in Y$. Define $Z_q=Z/(X\times\{y\})$ to be the identification of $X\times\{y\}$ with a point $q$. Then we say that $g:Z_q\longrightarrow Z_q$ is a \emph{ $q$-cut fiber injection} if $g|_{Z_q-\{q\}}$ is a fiber injection.

\end{definition}

\begin{theorem} If $g:B^S_p\longrightarrow B^S_p$ is a  $p$-cut fiber injection of a smooth fan $B^S_p$, then $g$ is not separated, continuum-wise semi-turbulent.
\end{theorem}
\begin{proof}
	If $B^S_p$ is a arc, then $B^S_p$ is homeomorphic to the blade space ${c}\times[0,1]$, $g$ is a homeomorphism and hence a fiber injection, so the result follows from Lemma 6.5.
	
	Suppose $B^S_p$ is not an arc and let $A, B$ be disjoint subcontinua of $B^S_p$. Then $p\not \in A$ or $p \not \in B$. Suppose the former. Then there is a $c_A\in S$ such that  $A \subset \nu_S(\{c_A\}\times(0,1])$ and the proof is the same as Lemma 6.5. 
\end{proof}

\vphantom{}

\begin{corollary}
	Let $h:Y\longrightarrow Y$ be a homeomorphism of a smooth fan $Y$. Then $h$ is not separated, continuum-wise semi-turbulent. Hence, smooth fans do not admit separated continuum-wise turbulent homeomorphisms.	

\end{corollary}

\begin{proof}
	This follows from the fact that homeomorphisms of smooth fans are $p$-cut fiber injections, where $p$ is the top of $Y$.
\end{proof}

We end with some questions:

\vphantom{}

\textbf{Question.} Does there exist a separated continuum-wise semi-turbulent homeomorphism on a non-smooth fan?

\vphantom{}

\textbf{Question.} If $F \subset X \times X$ is a closed relation on a compact metric space $X$ and $ent(F)>0$, is $(X_F^+, \sigma_F^+)$ either turbulent, reverse turbulent, or semi-turbulent?

\vphantom{}


\vphantom{}









\vphantom{}

\noindent  Judy  Kennedy\\
             Department of Mathematics,  Lamar University, 200 Lucas Building, P.O. Box 10047, Beaumont, Texas 77710 USA\\
{{kennedy9905@gmail.com}       }  \\

\noindent Christopher Mouron \\
              Department of Mathematics, Rhodes College, \\2000 North Parkway, Memphis, TN 38112 \\
                  {mouronc@rhodes.edu}           

\vphantom{}
     
\noindent Van Nall\\
             Department of Mathematics and Computer Science, University of Richmond, \\
               410 Westhampton Way, University of Richmond, VA 23172\\
               {vnall@richmond.edu}       				\-
				

\end{document}